\documentclass{proc-l}
\usepackage{amssymb}

\copyrightinfo{2000}{American Mathematical Society}

\newtheorem{Thm}{Theorem}[section]

\newtheorem{Lem}[Thm]{Lemma}
\newtheorem{Prop}[Thm]{Proposition}
\numberwithin{equation}{section}
\theoremstyle{definition}

\theoremstyle{remark}

\def\ldots{\mathinner{\ldotp\ldotp\ldotp}}
\def\ldots{\mathinner{\cdotp\cdotp\cdotp}}

\def \cal{\mathcal}

\begin{document}

\title{The range of operators on von Neumann algebras}

\author{Teresa Berm\'{u}dez}
\address{Departamento de An\'{a}lisis Matem\'{a}tico \\
Universidad de La Laguna \\
38271 La Laguna (Tenerife) \\
Canary Islands\\Spain.}

\email{tbermude@ull.es}

\author{N. J. Kalton}
\address{Department of Mathematics \\
University of Missouri-Columbia \\
Columbia, MO 65211\\ USA
}

\email{nigel@math.missouri.edu}

\thanks{The first author was supported by
DGICYT Grant PB 97-1489 (Spain); the second  was supported by NSF grant
DMS-9870027}

\begin{abstract}
We prove that for every bounded linear operator $T:X\to X$, where
$X$ is a non-reflexive quotient of a von Neumann algebra, the point
spectrum of $T^*$ is non-empty (i.e. for some $\lambda\in\mathbb C$ the
operator $\lambda I-T$ fails to have dense range.)
In particular, and as an application, we obtain that such a space
cannot support a topologically transitive operator.
\end{abstract}

\maketitle

\section{Introduction}

The results in this note are motivated by a question related to
hypercyclic operators.   In \cite{Godefroy-Shapiro} G. Godefroy and J.  Shapiro
suggest an extension of the notion of a hypercyclic operator to Banach
space which are not necessarily separable, via the notion of
topologically transitive operators  (see Section \ref{applications}
below).  Every Hilbert space supports a topologically transitive operator
(see the example due to J. Shapiro in Section \ref{applications}.)
Recently,
it has been shown by S. Ansari \cite{ansari} and L. Bernal \cite{bernal} that
every separable Banach space supports a hypercyclic operator, so it is
natural to ask whether every Banach space supports a topologically
transitive operator.

It is well-known that if $T$ is hypercyclic then the adjoint operator $T^*$ has
empty point
spectrum, $\sigma_p(T^*),$ \cite {herrero} and \cite{kitai}; this extends
to topologically
transitive operators (Proposition \ref{espectro-puntual}).  Thus we are
led to the question whether there exist complex Banach spaces so that
for every operator $T$ we have $\sigma_p(T^*)\neq\emptyset.$   Such
an example exists in the literature,
\cite{Shelah} and \cite{Shelah2}.  However, we show here that there are
much more natural
examples.  If $X$ is any von Neumann algebra (or even a non-reflexive
quotient of a von Neumann algebra), then any operator $T$ on $X$ has
$\sigma_p(T^*)\neq \emptyset.$  In particular this holds if
$X=\ell_{\infty}$ or $X=\mathcal L(\ell_2).$  We note hypercyclicity
with respect to the strong-operator topology on $\mathcal L(\ell_2)$
has been considered in \cite{chan} and \cite{montes}.

Our main result is rather stronger in that we show that if $X$ is
a non-reflexive quotient of a von Neumann algebra, then for any operator
$T$ we have that the quotient space $X/\overline{{\mathcal
R}(\lambda-T)}$ contains a copy of $\ell_{\infty}$ and is in particular
non-separable.

Let us point out by way of further motivation that  any
operator $T$ on $\ell_1$ satisfies $\sigma_p(T^{**})\neq\emptyset$, since
if $\lambda$ is in the approximate point spectrum of $T$ then it is in
the point spectrum of $T^{**}$ by an argument depending on the Schur
property of $\ell_1$
(this was shown to us by M. Gonz\'{a}lez).  This is suggestive of the main
result in the case $X=\ell_{\infty}.$

Our arguments depend on two Banach space concepts, which we now
introduce.  A projection $P$ on a Banach space $X$ is an L-projection if
$\|x\|=\|Px\|+\|x-Px\|$ for any $x\in X.$ A Banach space
$X$ is said to be
{\it
L-embedded} if there is an L-projection of $X^{**}$ onto $X$ i.e. if
there is a projection
$\Pi:X^{**}\to X$ so that we have:
$$ \|x^{**}\|=\|x^{**}-\Pi x^{**}\|+\|\Pi x^{**}\| \qquad  \mbox{ for }  x^{**}\in
X^{**}.$$   For the basic facts on
L-embedded spaces we refer to
\cite{HWW} Chapter IV.  A Banach space $X$ is called a {\it Grothendieck
space} if every bounded operator $T:X\to Y$ with separable range is
weakly compact. This is equivalent to requiring that if $\{x_n^*\}_{n\in \mathbb N} $ is a
weak$^*$-null sequence in $X^*$, then it is also weakly null.  Any von
Neumann algebra is a Grothendieck space \cite{P} and its dual is
L-embedded \cite{Takesaki}, \cite{HWW}.  We also recall that a Banach
space $X$ is called an {\it Asplund space} if every separable subspace
has separable dual (this is equivalent to the original definition,
\cite{DGZ} Theorem 5.7, p.29).

Most of our notation is standard.
We will use $B_X$ to denote the closed unit ball of a Banach space
$X$.  If $F$ is a subset of $X$ then $\langle F\rangle$ denotes its
linear span.

We would like to thank M. Gonz\'{a}lez, J. Shapiro and D. Werner for helpful
comments.

\section{Main results}

We use repeatedly the following principle:

\begin{Lem}\label{lema}\cite[II.E.15]{Wojtaszczyk}
 Let $X$ be a Banach space and suppose
$\{C_k\}_{k=1}^n$ is a
 finite set of  convex sets.  Suppose $D_k$ is the
weak$^*$-closure of $C_k$ in $X^{**}.$  If
$\cap_{k=1}^nD_k\neq \emptyset$ then for any $\epsilon>0$ there exists
$x\in C_1$ with
$d(x,C_k)<\epsilon$ for $k=2,3\ldots,n.$
\end{Lem}

We will also need the following well-known variant of the
Hahn-Banach Theorem.

\begin{Lem}\label{HB} Let $X$ be a Banach space and suppose $F$ is a
finite-dimensional subspace of $X^*$.  If $\psi$ is a linear functional
on $F$ with $\|\psi\|<1$ then there exists $x\in X$ with $\|x\|<1$ and
$x^*(x)=\psi(x^*)$ for $x^*\in F.$ \end{Lem}

\begin{proof} This can be proved directly or from Lemma \ref{lema}.  Let
$C_1=\{x\in X:\ x^*(x)=\psi(x^*)\ \forall x^*\in F\}$ and $C_2=\{x\in X:
\
\|x\|\le
\|\psi\|\}.$  Then, by the Hahn-Banach Theorem, the weak$^*$-closure
$D_1$ of
$C_1$ is the set
$\{x^{**}\in X^{**}:\ x^{**}(x^*)=\psi(x^*)\ \forall x^*\in F\}$.  By an
application of the Hahn-Banach Theorem and Goldstine's Theorem
(\cite{megginson} Theorem 2.6.26, p. 232) $D_1$
meets the weak$^*$-closure of $C_2$ so that we can apply Lemma
\ref{lema}.\end{proof}

\begin{Prop}\label{oneone}
Suppose $T:X\to Y$ is a bounded linear
operator.  Then the following properties are equivalent:\newline
(1) ${\mathcal N}(T^{**})=\{0\}.$\newline
(2) If $\{x_n\}_{n\in \mathbb N}\subset X$ is a bounded sequence such
that
$\displaystyle\lim_{n\to\infty}\|Tx_n\|
=0$ then $\displaystyle\lim_{n\to\infty}x_n=0$ weakly.
\end{Prop}

\begin{proof} (1) implies (2).  Clearly $0$ is the only
weak$^*$-cluster point of $\{x_n\}_{n\in \mathbb N}$ in $X^{**}$ and so
$\displaystyle\lim_{n\to \infty} x_n=0$ weakly.

 (2) implies (1).  Assume for some $x^{**}\neq 0$ with $\| x^{**}\|=1$, we
have
$T^{**}x^{**}=0.$  Pick $x^*\in X^{*}$ with
$x^{**}(x^*)=1.$
 Then for each $n$ the sets $C_1=\{x:
\|x\|\le 1\}$, $C_2=\{x: x^*(x)\ge 1\}$
and $C_3=\{x:\|Tx\|\le n^{-1}\}$ satisfy the
conditions of  Lemma \ref{lema}, so we pick $\{x_n\}_{n\in
\mathbb N}\subset X$ with
$\|Tx_n\|\le n^{-1}$,
 $\|x_n\|\le 2$ and $x^*(x_n)\ge \frac12,$
contradicting (2).\end{proof}

Now if $T:X\to Y$ is a bounded linear  operator we denote by
$\hat T$ the induced operator $\hat T:X^{**}/X\to
Y^{**}/Y.$

\begin{Prop}\label{lowerbound} Suppose $T:X\to Y$ is a
bounded operator.  Then the following are equivalent:
\newline
(1) There exists a sequence $\{\xi_n\}_{n\in \mathbb N}\subset  X^{**}/X$ such
that $\|\xi_n\|=1$ and  \linebreak $\displaystyle \lim_{n\to \infty}\|\hat
T\xi_n\|=0.$\newline
(2) There exists a bounded sequence $\{x^{**}_n\}_{n\in
\mathbb N}\subset
X^{**}$ such that $d(x^{**}_n,X)=1$ and $\displaystyle\lim_{n\to \infty}
\|T^{**}x_n^{**}\|=0.$\end{Prop}

\begin{proof}  We only need prove (1) implies (2).
Pick $w_n^{**}\in\xi_n$ with $\|w_n^{**}\|\le 2.$  Let
$\epsilon_n:=\|\hat T\xi_n\|+\frac1n.$  Then there exists
$u_n\in X$ with $\|T^{**}w_n^{**}-u_n\|
<\epsilon_n.$  We now argue that since
$T^{**}w_n^{**}$ is in the  weak$^*$-closure of both
$u_n+ \epsilon_nB_X$ and $2T(B_X)$, then  there exists $v_n\in X$ with
$\|v_n\|\le 2$ and $\|Tv_n-u_n\|\le 2\epsilon_n.$  Thus
$\|T^{**}(w_n^{**}-v_n)\|\le 3\epsilon_n.$  Letting
$x_n^{**}:=w_n^{**}-v_n$ we are done.\end{proof}

\begin{Thm}\label{lowerbound2} Suppose $X$ is a
subspace of an L-embedded Banach space $V$, and $Y$ is any Banach space.
 Suppose $T:X\to Y$ is a
bounded linear  operator such that ${\mathcal N}(T^{**})\subset X.$
 Then there exists $\delta>0$ so that for all $\xi\in
X^{**}/X$ we have $\|\hat T\xi\|\ge
\delta\|\xi\|.$\end{Thm}

\begin{proof}
We start by proving the Theorem in the special case when ${\mathcal
N}(T^{**})=\{0\}.$

  Suppose the conclusion is false. Using
Proposition \ref{lowerbound} we produce a bounded
sequence $\{x_n^{**}\}_{n\in \mathbb N}\subset  X^{**}$ with
$d(x_n^{**},X)=1$ but $\lim_{n\to \infty } \|T^{**}x_n^{**}\|=0.$
We can regard  $X^{**}$ as a subspace of
$V^{**}$.  Now let $\delta_n=d(x_n^{**},V).$  For fixed $n$, if
$\rho>\delta_n$, then
$x_n^{**}$ is in the weak$^*$-closure of both $X$ and $v+\rho B_V$ for
some $v\in V.$  Hence there is $y\in v+\rho B_V$ such that
 $d(y,X) \le \rho$ by an application of Lemma
\ref{lema} and so $d(x_n^{**},X)\le 2\rho.$  We conclude that
$\delta_n\ge \frac12$ for each $n\in\mathbb N.$

Let us denote by $\Pi$ the $L$-projection of $V^{**}$ onto $V,$ and let
$V_s=\text{ker }\Pi.$  Let
$v_n:=\Pi
x_n^{**}$ and $v_n^{**}:=x^{**}_n-v_n.$  Then $v_n^{**}\in V_s$ and
$\|v_n^{**}\|=\delta_n\ge \frac12.$ Let $a:=\sup _{n\in \mathbb N}
\|x^{**}_n\|$ and $\eta_n:=\|T^{**}x_n^{**}\|+\frac1n$.

We shall define inductively a sequence
$\{x_n\}_{n\in \mathbb N}$ in $X,$ and a sequence $\{x_n^*\}_{n\in \mathbb N}$ in
$X^*$ such that:
\begin{equation}\label{one}
\|x_n\|\le a \qquad n\in\mathbb N
\end{equation}
\begin{equation}\label{onea}
\|Tx_n\|<\eta_n
\end{equation}
\begin{equation}\label{two}
 \|x^*_n\|<1
\qquad n\in\mathbb N
\end{equation}
\begin{equation} \label{three}
|x^*_n(x_k)|\ge \frac18
  \quad 1\le k\le n .
  \end{equation}

Let us suppose $n\in\mathbb N$ and
$\{x_k\}_{k<n},$ and
$\{x^*_k\}_{k<n}$ have been determined and
satisfy
(\ref{one}), (\ref{onea}), (\ref{two}) and (\ref{three});
 if $n=1$ these sets
are empty of course. We shall determine
$x_n$ and $x_n^*$.

Let $ F:=\langle \{x_1,\ldots,x_{n-1},v_n\}\rangle$ and
$G:=\langle \{x_1,\ldots,x_{n-1},v_n,v_n^{**}\}\rangle.$
 If $n>1$  we define
$\psi=\psi_n\in F^*$ by taking
$\psi$ to be a norm-preserving
extension of
$x^*_{n-1}|_{F\cap X};$
if
$n=1$ we simply let
$\psi=0.$
Then $\|\psi\|<1.$  Let $\psi(v_n)=re^{i\theta}$ where
$0\le \theta<2\pi$ and $r\ge 0.$
We next define
$\varphi\in G^*$ to the
extension of
$\psi$ such that
$\varphi(v_n^{**})=\frac14e^{i\theta}.$
We claim that $\|\varphi\|<1.$    In fact if
$u^{**}\in G$ then we
can write $u^{**}=u+ \mu v_n^{**}$ where $\mu\in\mathbb C$ and $u\in F$.
Then
\begin{align*} |\varphi(u^{**})|&\le
|\psi(u)|+\frac14|\mu|\\
&\le
\|\psi\|\|u\|+\frac12 |\mu|\|v_n^{**}\| \\ &\le
\max(\frac12,\|\psi\|)\|u^{**}\|<\| u^{**}\|.\end{align*}
Now by Lemma \ref{HB} we can define
$v^*\in V^*$
 with
$\|v^*\|<1$ and
$u^{**}(v^*)=\varphi(u^{**})$
for
$u^{**}\in G.$  Let
$x^*_{n}$ be the restriction of
$v^*$ to $X$.

Now consider the sets $C_1=\{x:\|x\|\le a\},\ C_2=\{x:\|Tx\|\le
\|T^{**}x_n^{**}\|\}$ and $C_{3}=\{x:\
x^*_{n}(x)=x^{**}_n(x^*_{n})
\}.$   Clearly $x^{**}_n$ belongs to the weak$^*$-closure of each set.  By
Lemma
\ref{lema} we can find $x_n\in C_1$ with $\|Tx_n\|<\eta_n,$ and so that
$$  |x^*_{n}(x_n)|> |x^{**}_n(x^*_n)|-\frac18.$$
It is now clear that (\ref{one}), (\ref{onea}) and (\ref{two}) hold.
For (\ref{three}) note that if $k<n$ we have $x^*_n(x_k)=x^*_{n-1}(x_k)$
while
$$ |x_n^*(x_n)|\ge |x_n^{**}(x_n^*)|-\frac18 = (\frac14+r)-\frac18\ge
\frac18.$$

Now the proof is completed (for the special case ${\mathcal
N}(T^{**})=\{0\}$) by observing if $x^*$ is any weak$^*$-cluster point of
the sequence $\{x_n^*\}_{n\in \mathbb N}$ then $|x^*(x_n)|\ge \frac18$ for all $n.$  Since
$\lim_{n\to\infty}\|Tx_n\|=0$ this contradicts Proposition \ref{oneone},
 since $x_n$ does not converge weakly.

To treat the general case suppose $R={\mathcal N}(T^{**})={\mathcal
N}(T).$ Then
$R$ is reflexive.  Consider the induced map $T_0:X/R\to Y$; clearly
${\mathcal N}(T_0^{**})=\{0\}.$  We next note that $X/R$ embeds into
$V/R$ and $V/R$ is L-embedded \cite {HWW} p.160.  Hence $\hat T_0$
satisfies a lower bound
on $Z=(X/R)^{**}/(X/R).$  However it is easily seen that  $Z$ coincides
with $X^{**}/X$ and $\hat T_0=\hat T.$\end{proof}

We next need some facts about Grothendieck spaces:

\begin{Prop}\label{Kernels}
Suppose $Y$ is a
Grothendieck space and that $T:X\to Y$ is a bounded
linear operator such that $T^{*}$ is one-one.  Then $T^{***}$
is one-one.
\end{Prop}

\begin{proof}  Suppose $\{y_n^*\}_{n\in \mathbb N}\subset  Y^*$ is a
bounded
 sequence  such that
$\lim_{n\to \infty }\|T^*y_n^*\|=0.$  Let $y^*$ be any
weak$^*$-cluster point of $\{y_n^*\}_{n\in \mathbb N}$.  Then $T^*y^*=0$
so that $y^*=0.$  Therefore $\displaystyle \lim_{n\to\infty} y_n^*=0$ weak$^*$.
But since $Y$ is a Grothendieck space this implies
$\displaystyle \lim_{n\to \infty } y_n^*=0$ weakly and we can apply Proposition
\ref{oneone}.
\end{proof}

\begin{Prop}\label{Grothendieck}  Suppose $X$ is a Grothendieck space and
$Y$ is a subspace of $X$ so that $X/Y$ is reflexive.  Then $Y$ is a
Grothendieck space.\end{Prop}

\begin{proof} Suppose $T:Y\to c_0$ is any bounded operator.  Then we may
find a Banach space $E\supset c_0$ with $E/c_0\cong X/Y$ and an extension
$\tilde T: X\to E.$  We claim $E$ is an Asplund space.  Indeed if $F$ is
a separable subspace of $E$ then let $F'$ be the closure of $c_0+F$ which
is also separable.  Then $F'/c_0$ is separable and reflexive so that
since $c_0^*\cong \ell_1$ is separable, $F'$ has separable dual.  Now it
follows from a deep result of Hagler and Johnson \cite{haglerjohnson}
(see also \cite{diestel})
that $B_{E^*}$ is weak$^*$-sequentially compact.  Hence if $(e_n^*)$ is
any sequence in $B_{E^*}$ there is a subsequence $(f_n^*)$ so that
$\tilde T^*f_n^*$ is weak$^*$ and hence weakly convergent in $X^*.$  Thus
$\tilde T$ is weakly compact by Gantmacher's theorem (see
\cite{megginson}  Theorem 3.5.13 p. 343) and in particular
$T$ is weakly compact.\end{proof}

\begin{Thm}\label{reflexive} Suppose $X$ and $Y$ are Banach spaces and
$Y$ is a Grothendieck space.  Suppose $T:X\to Y$ is a bounded operator
such that $Y/\overline{{\mathcal R}(T)}$ is reflexive.  Then ${\mathcal
N}(T^{***})\subset Y^*.$ \end{Thm}

\begin{proof} Let $Y_0=\overline{{\mathcal R}(T)}.$  Then by
Proposition \ref{Grothendieck} $Y_0$ is also a Grothendieck space.  We
write $T=JT_0$ where $J:Y_0\to Y$ is the inclusion map and $T_0:X\to Y_0$.   Clearly
$(Y/Y_0)^*\cong
{\mathcal N}(T^*)$ is reflexive. We observe that $T_0^*$ is one-one and
by Proposition \ref{Kernels} we obtain that $T_0^{***}$ is also one-one.
Now, since $Y/Y_0$ is reflexive, this implies  ${\mathcal N}(T^{***})=
{\mathcal
N}(J^{***})={\mathcal N}(J^*)\subset Y^*$ as required.\end{proof}

\begin{Thm}\label{main} Let $X$ be a non-reflexive complex Banach space
which is a Grothendieck space
such that $X^*$ is isometric to a subspace of an $L$-embedded space.
Suppose $T:X\to X$ is a bounded linear operator.  Then
there exists $\lambda\in\mathbb C$ so that $X/\overline{{\mathcal
R(\lambda -T)}}$ is non-reflexive (and hence non-separable).
  In particular the point spectrum $\sigma_p(T^*)$ is
non-empty.\end{Thm}

 \begin{proof} Let $S=T^*$.  Then since $X$ is non-reflexive, the operator
$\hat S$ has non-empty spectrum and furthermore for any
$\lambda$ in the boundary
$\partial \sigma(\hat S)$ there  is a sequence $\xi_n\in
X^{***}/X^*$ with $\|\xi_n\|=1$ so that $\lim_{n\to\infty}\|
(\lambda -\hat S)\xi_n\|=0.$
This implies that for $\lambda\in\partial \sigma (\hat S)$ we have
${\mathcal N}((\lambda-S)^{**})$ is not contained in $X^*$ by Theorem
\ref{lowerbound2}.  Then we apply Theorem \ref{reflexive} and deduce that
$X/\overline{{\mathcal R}(\lambda -T)}$ is non-reflexive.
By Proposition \ref{Kernels} we have that  $\lambda \in  \sigma _p(T^*)$.\end{proof}

Our main example for Theorem \ref{main} is when $X$ is a von Neumann
algebra.  The fact that von Neumann algebras have the Grothendieck
property is a recent result of Pfitzner \cite{P}.  In fact slightly more
follows from Pfitzner's work.

\begin{Prop}\label{Pfitzner}  Let $A$ be a von Neumann algebra and
suppose $T:A\to Y$ fails to be weakly compact, then there is a closed
subspace $E$ of $A$ such $T|_E$ is an isomorphism and $E$ is isomorphic
to $\ell_{\infty}.$\end{Prop}

\begin{proof} Suppose $T$ fails to be an isomorphism on any subspace
isomorphic to $\ell_{\infty}.$ Let
$A_0$ be any maximal Abelian subalgebra of
$A.$ Then
it follows from classical results of Rosenthal \cite{R}
 that $T$ is weakly compact on $A_0$; by Pfitzner's Theorem
\cite{P} Theorem 1 (see also Corollary 10), $T$ is weakly
compact.\end{proof}

\begin{Thm}\label{main2} Let $X$ be a non-reflexive quotient  of a von
Neumann algebra, and let $T:X\to X$ be any bounded linear operator.  Then
there exists $\lambda\in\mathbb C$ so that $X/\overline{{\mathcal
R(\lambda -T)}}$ contains an isomorphic  copy of $\ell_{\infty}$ and
hence ${\mathcal N}(\lambda -T^*)$ contains an isomorphic copy of
$\ell_{\infty}^*.$  In particular the point spectrum $\sigma_p(T^*)$ is
non-empty.\end{Thm}

\begin{proof}
The dual of any $C^*$-algebra is $L$-embedded (\cite{Takesaki},
\cite{HWW}) and so it follows from the work of Pfitzner \cite{P} that
$X$ satisfies the hypotheses of Theorem \ref{main}.  Proposition
\ref{Pfitzner} implies that if $X/\overline{{\mathcal
R(\lambda -T)}}$  is non-reflexive then it contains a
complemented isomorphic
copy of $\ell_\infty$. Since $\left(X/\overline{{\mathcal
R(\lambda -T)}}\right) ^*\cong {\mathcal N}(\lambda -T^*)$, there
exists in ${\mathcal N}(\lambda -T^*)$ an isomorphic copy of
$(\ell_\infty)^*$.
\end{proof}

\section{Applications to hypercyclic operators}\label{applications}

A bounded linear operator $T$ on a complex Banach space $X$ is
called {\em hypercyclic} if there is a vector $x\in X $
 (called {\em hypercyclic vector for $T$}) such that
 $\{ T^nx\;\;:\;\; n\in \mathbb N\}$  is dense on $X$.  This concept is
related to the problem of the existence of proper closed invariant
subsets for
a bounded linear operator.  It is an open problem whether every bounded
linear operator on a Hilbert space has a proper closed invariant subset,
or equivalently if every operator has a non-zero vector which is not
hypercyclic.  We refer to \cite{Grosse} for an excellent survey.

We note that a non-separable Banach space cannot support a hypercyclic
vector.
An approach to obtain something similar to  hypercyclicity  in a non-separable
Hilbert and Banach  spaces was given by K. Chan \cite{chan} and
A. Montes and C. Romero \cite{montes}, respectively. In fact, they give  certain
``hypercyclicity" results  in ${\mathcal L}(X)$ where $X$ is a
separable Banach space, using the strong operator topology in place of
the standard uniform norm topology.

It is however possible to extend the notion of hypercyclic operators to
nonseparable Banach spaces in a natural way using result of
\cite{Godefroy-Shapiro}. Let us say that an operator $T$ on an arbitrary
Banach space is {\em topologically transitive} if for every pair $U,\;V$ of non-void open subsets of $X$, there
exists a positive integer $n$ such that $T^n(U) \cap V\neq
\emptyset. $
 In Theorem 1.2 of \cite{Godefroy-Shapiro} it is proved
 if $X$ is a
separable  Banach space, then $T$ is hypercyclic if and only if   $T$ is
topologically transitive.

\ \par

The following Proposition is immediate.

\begin{Prop}\label{Toptrans}  A bounded linear operator $T$ is
topologically transitive if and only if every proper closed invariant
subset has empty interior.\end{Prop}

An  argument similar to the result due to J. B\'{e}s and A. Peris
\cite{Bes-Peris} provides  a
sufficient condition  for topological transitivity.

%The following Proposition, which provides a criterion for topological
%transitivity is proved in a similar way to a result
%due to B\'{e}s and Peris
%\cite{Bes-Peris}.

\begin{Prop}\label{hc}
{\em (Topologically transitive Criterion)}
Let $T$  be a bounded linear operator  on a complex Banach space
$X$ (not necessarily separable). Suppose that there exists a strictly
increasing sequence of
 positive integers $\{ n_k\}_{k\in \mathbb N}\subset  \mathbb N$  for which there
 are:
\begin{enumerate}
  \item A dense subset $X_0\subset X$ such that $ T^{n_k}x\to 0$ for every
  $x\in X_0$.
  \item A dense subset $Y_0\subset X$ and  a sequence of mappings
  $S_{k}:Y_0 \to X$ such that
\begin{enumerate}
  \item $ S_{k}y\to 0$ for every
  $y\in Y_0$.
  \item $T^{n_k}S_{k}y\to y$ for every $y\in Y_0$.
\end{enumerate}
 \end{enumerate}
Then $T$ is topologically transitive.
\end{Prop}

{\it Example.}  The following example was
suggested by J. Shapiro. Let us use Proposition
\ref{hc} to show that there is a topologically transitive operator on any
Hilbert space. If $H$ is separable, then the result is clear by
S. Ansari \cite{ansari}
and
L.  Bernal \cite{bernal}.
 If $H$ is a non-separable Hilbert space we write $H=\ell_2(X)$
where $X$ is a Hilbert space of the same density character.
Define $T$ as twice the backward shift on $\ell_2(X)$, that is
$$
T(x_1, x_2, \dots ):=2(x_2, x_3, \dots ).
$$
Using Proposition \ref{hc}, we have  that $T$ is topologically
transitive
 taking  $n_k=k,$
$$
X_0:=\{ \mbox{ finitely non-zero sequences in } \ell_2(X)\}
$$
$$
Y_0:=\ell_2(X),
$$
$$
S(x_1, x_2, \dots ):=\frac{1}{2} (0, x_1, x_2, \dots )
$$
and $S_{k}:=S^{k}$.

Clearly this example can be modified to replace $H$ by any space
$\ell_p(I)$ where $1\le p<\infty.$

It has been shown  by S. Ansari \cite{ansari}
and
L.  Bernal \cite{bernal} that any separable complex
Banach space supports a hypercyclic operator. Recently, J. Bonet and A. Peris
gave a version for  ${\cal F}$-spaces
\cite{bonet-peris}.
This suggests the corresponding problem of determining whether every
complex Banach
space supports a topologically transitive operator.

This question has a negative answer.
   In order to see this, we need to give a spectral
property of  topologically transitive operators, which is well-known in
the case of hypercyclic operators \cite{kitai} and \cite{herrero}.

\begin{Prop}\label{espectro-puntual}
Let $T$ be a bounded linear operator on a complex Banach space. If
$T$ is topologically transitive, then  $\sigma_p(T^*)$ is empty.
\end{Prop}

\begin{proof} If $\lambda\in\sigma_p(T^*)$ and $x^*$ is a corresponding
eigenvector then one of the sets $\{x:|x^*(x)|\ge 1\}$ or
$\{x:|x^*(x)|\le 1\}$ is an invariant set with non-empty interior. Then
use Proposition \ref{Toptrans}.\end{proof}

As pointed out in the introduction the examples of \cite{Shelah} and
\cite{Shelah2} of non-separable spaces such that every bounded operator is
a perturbation of a multiple of the identity by an operator with
separable range give examples of spaces which support no topological
transitive operators.  However, the following theorem shows that
$\ell_{\infty}$ and
$\mathcal L(\ell_2)$ are more natural examples, where $\mathcal L(\ell_2)$ denotes
the algebra of all bounded linear operators on $\ell_2.$

\begin{Thm}\label{appli}  Let $X$ be a non-reflexive quotient of a von
Neumann algebra.  Then $X$ does not support a topologically transitive
operator.  In particular $\mathcal L(\ell_2)$ and $\ell_{\infty}$ do not
support a topologically transitive operator.
\end{Thm}

\begin{proof}  Just apply Theorem \ref{main2} and Proposition
\ref{espectro-puntual}.\end{proof}

\ \par

We conclude with a remark on ultrapowers.
We recall some concepts about ultrapowers of Banach spaces
and operators. See \cite{Heinrich}
for more information.
We fix a non-trivial ultrafilter $\mathcal U$ on the set $\mathbb N$ of all
positive integers.
For every Banach space $X$, we
consider the Banach space $\ell_{\infty}(X)$
of all bounded sequences $(x_n)$ in $X$, endowed with
the norm $\|(x_n)\|_{\infty}:= \sup\{\|x_n\|: n\in \mathbb N\}$.
Let  $N_{\mathcal U}(X)$ be  the closed subspace of all sequences
$(x_i)\in\ell_{\infty}(X)$ which converge to $0$ following $\mathcal U$.
The {\em ultrapower of $X$ following $\mathcal U$} is defined as the
quotient
$$
X_{\mathcal U}:=\frac{\ell_{\infty}(X)}{N_{\mathcal U}(X)}.
$$
The element of $X_{\mathcal U}$ including the sequence
$(x_i)\in \ell_{\infty}(X)$ as a representative is
denoted by $[x_i]$. Its norm in $X_{\mathcal U}$ is given by
$$
  \bigl\|[x_n]\bigr\|=\lim_{\mathcal U}{\| x_n\|}.
$$
The constant sequences generate a subspace of $X_{\mathcal U}$
isometric to $X$. So we can consider the space $X$ embedded in
$X_{\mathcal U}$. Moreover, every operator $T\in L(X, Y)$ admits an
extension $T_{\mathcal U}\in L(X_{\mathcal U},Y_{\mathcal U})$, defined by
$$
T_{\mathcal U}([x_n]) := [T x_n],\,\,\, [x_n]\in X_{\mathcal U}.
$$

An easy argument with ultrapowers gives that any ultrapower cannot
be a topologically transitive operator. This fact  can be obtained  by
 the following  easy argument.
\begin{Prop}
Let $\mathcal U$ be an ultrafilter, $X$ a
complex Banach space and $T$ any
bounded linear operator on $X$. Then $T_\mathcal U$ is not topologically
transitive.
\end{Prop}

\begin{proof}  We note that any $\lambda\in\partial \sigma(T)$ is
in the approximate point spectrum of $T^*$, i.e. there exists a sequence
$\{x_n^*\}_{n\in \mathbb N}$ in $X^*$ with $\|x_n^*\|=1$ and
$\lim_{n\to\infty}\| \lambda x_n^*-T^*x_n^*\|=0.$  Now let $\xi^*\in X_{\mathcal
U}^*$ be defined by $\xi^*([x_n])=\displaystyle \lim_{\mathcal U}x_n^*(x_n).$  Then
$\xi^*\in\sigma_p(T_{\mathcal U}^*)\neq \emptyset,$ so we can apply
Proposition
\ref{espectro-puntual}.\end{proof}

\ \par

We conclude with the following open question: {\em Is there any
characterization of non-separable Banach spaces which support a
topologically transitive operator?}

 \end{document}